\documentclass[11pt]{elsarticle}

    \usepackage{color}  
\usepackage{lineno,hyperref}
\usepackage{amsmath,amsthm,amssymb,amsfonts}
\usepackage{mathrsfs}
\usepackage{subfig,graphicx}
\modulolinenumbers[5]
\allowdisplaybreaks[3]

\journal{International Journal of Control}

\newtheorem{theorem}{Theorem}[section]
 
  \newtheorem{lemma}[theorem]{Lemma}
   
    \theoremstyle{definition}

\theoremstyle{remark}
 \numberwithin{equation}{section}

\begin{document}

\begin{frontmatter}

\title{Unknown input observer design for a class of coupled 
wave PDE systems\footnote{\textcolor{red}{\underline{\textbf{To cite this article:}}} \textcolor{blue}{\textbf{Ghaderi,~N., \& Jacob, B. (2026) Unknown input observer design for a class of coupled  wave PDE systems, International Journal of Control, DOI:10.1080/00207179.2025.2612055}}\\ \underline{To link to this article:} \textbf{https://doi.org/10.1080/00207179.2025.2612055}}} 

\author[mymainaddress,secondaddress]{Najmeh Ghaderi\corref{mycorrespondingauthor}}
\cortext[mycorrespondingauthor]{Corresponding author}
\ead{najmeh.ghaderi@uni-saarland.de}

\author[mymainaddress]{Birgit Jacob}
\ead{bjacob@uni-wuppertal.de}

\address[mymainaddress]{School of Mathematics and Natural Sciences, University of Wuppertal, Gaußstr. 20, 42097 Wuppertal, Germany}
\address[secondaddress]{Chair of Systems Theory and Control Engineering, Saarland University, Saarbrücken, Germany}

\begin{abstract}
This paper deals with the problem of designing unknown input observers for a class of coupled semilinear wave partial differential equations (PDE) systems. A state observer is designed to estimate the uncertain coupled wave PDE systems. Then, the analysis of the asymptotic stability and $H_{\infty}$ performance for the observer design of coupled wave PDE systems is investigated. Some sufficient conditions of asymptotic stability for the observer error system with disturbance attenuation level are derived via matrix inequalities based on the Lyapunov stability theory. Finally, a numerical simulation is presented to demonstrate the effectiveness of the obtained result. 
\end{abstract}


\end{frontmatter}

\section{Introduction}
The behavior of various complex physical processes is described through partial differential equations (PDEs), which model the development of spatially distributed phenomena \cite{cristofaro2}. Such examples can be seen in hydraulic networks, tubular chemical reactors, etc.~\cite{dos,luyben}. Among the wide range of PDE equations, the wave equation is one of the most important ones in applications like mechanics, as it describes the movement in both strings and wires, as well as the fluid dynamics of surfaces, e.g., water waves \cite{rapp}. State estimation and stabilisation of PDE systems are the crucial issues in control theory. Therefore, the stabilization and estimation of PDE systems have gained increased focus during the last two decades, e.g., \cite{curtain,ghaderi,jacob1,jacob2,krstic,tucsnak}.

State estimation in the presence of unknown inputs is also a critical challenge in control theory due to the inevitable uncertainties in systems, including system parameter variation, external disturbances, and fault signals, which can all be regarded as unknown inputs impacting the system \cite{zhu}. Indeed, in real-world applications, some state variables might not be directly measurable or could be influenced by external disturbances or inputs that are not fully known. Therefore, an unknown input observer (UIO) is used to estimate the system state in the presence of unknown inputs \cite{commault,hammouri,luenberger,zhang}. In addition, over the past few decades, various methods have emerged for designing observers of ordinary differential equations (ODEs) to tackle uncertainties in systems \cite{chen,commault,corless,delshad,hammouri}. Furthermore, $H_{\infty}$ observer design plays a crucial role in estimating the unmeasured states of a system while ensuring robustness against disturbances and uncertainties. For instance, the work by \cite{du,gahinet,pertew} has explored $H_{\infty}$ observer design strategies for different system types of ODEs, contributing to the understanding of observer robustness in the face of uncertainties. The problem of $H_{\infty}$ and UIOs for distributed parameter systems has been recently studied in some works, for example \cite{demetriou,pisano} have extended the finite-dimensional UIO design approaches to a class of parabolic PDE systems, and observer design for wave equations with unknown periodic input has been considered in \cite{chauvin}. Furthermore, \cite{cristofaro1} addressed UIO design for nonlinear systems with full and partial information, utilizing the theory of characteristics to solve PDEs. \cite{li} developed an observer for Lipschitz nonlinear parabolic PDE systems, and \cite{cristofaro2} presents the design of a UIO for coupled first-order hyperbolic PDE/ODE linear systems with unknown boundary inputs. However, there is still much work to be done in this area. 
 
Inspired by \cite{li}, this paper aims to explore the observer design for a class of coupled wave systems with unknown inputs. We are going to establish the sufficient conditions necessary for ensuring the asymptotic stability of this design approach. The key contributions of this study are outlined as (i) We introduce an observer design technique for a class of uncertain wave PDE systems. The study examines the conditions for the existence and asymptotic stability of the observer when $d(t)=0$. (ii) The design of an $H_{\infty}$ observer is analyzed for uncertain semilinear coupled wave PDE systems under $d(t)=0$ and $d(t) \neq 0$, with $ \| e(x, t) \|_{\mathcal{L}^{2}} < \mu \| d(t) \|_{\mathcal{L}^{2}}$. The asymptotic stability is investigated, and design conditions are expressed through matrix inequalities. 
\section{Problem statement}
We consider a system of $n$ semilinear one dimensional wave equations with unknown inputs (disturbances)
\begin{align}\label{main}
\begin{split}
    w_{tt}(x,t) & =A w_{xx}(x,t)- B w_{t}(x,t)-D w(x,t) \\
    & \quad +f(w(x,t))+G u(t)+F d(t), \\
    w_{x}(0,t)& = C_{1} w_{t}(0,t),\\
    w_{x}(1,t)& =0,\\
    y(t)& =C \int_{0}^{1} w(x,t) \, dx+K u(t)+H  d(t),
    \end{split}
\end{align}
where $(x,t) \in (0,1) \times (0, \infty)$, $ w(x,t)= (w_1(x, t), w_2(x, t), \cdots, w_{n}(x, t))^\top \in [\mathcal{L}^{2}\big( (0,1); \mathbb R \big)]^{n}$
stands for the system state vector,
  $u(t) \in \mathbb R^{p}$ is the known input, $d(t) \in \mathbb R^{d}$ is the unknown input, and $y(t)  \in \mathbb R^{q}$ is a measured output. The matrices $A \in \mathbb R^{n \times n}$, $B \in \mathbb R^{n \times n}$, $C \in \mathbb R^{q \times n}$, $D \in \mathbb R^{n \times n}$, $C_{1} \in \mathbb R^{n \times n}$, $G \in \mathbb R^{n \times p}$, $F \in \mathbb R^{n \times d} $, $K \in \mathbb R^{q \times p} $ and $H \in \mathbb R^{q \times d} $ are known constant matrices, where $A$, $B$, $D$, and $C_{1}$ are positive definite matrices. The initial values are $w_{0}(x)=w(x,0)$ and $w_{1}(x)=w_{t}(x,0)$,  with $w_0(x)=( w_{10}(x),  w_{20}(x), \cdots, w_{n0}(x))^\top$, and $w_1(x)=(w_{11}(x),  w_{21}(x), \cdots, w_{n1}(x))^\top$. We assume that the nonlinear function $f$ is globally  Lipschitz continuous on $[\mathcal{L}^{2}\big( (0,1); \mathbb R \big)]^{n}$, that is, there exists a constant $\gamma>0$ such that
  $$   \| f(v)-f(w) \|_{\mathcal{L}^{2}} \leq \gamma \| v - w \|_{\mathcal{L}^{2}}, ~~ \forall v, w \in [\mathcal{L}^{2} \big((0,1);\mathbb R\big)]^{n}.$$
  
 The following lemma will be used later.
\begin{lemma} (\cite{cao})
Let the set of $n$ by $n$ positive definite symmetric matrices be denoted by $\mathcal{S}_{n}^{+}$. Then for  $x, y \in \mathbb R^{n}$ and $S \in \mathcal{S}_{n}^{+}$, we have
\begin{align}
&  2 x^\top  y \leq x^\top  S^{-1} x + y^\top  S y,\label{young1}\\
& - 2 x^\top  y \leq  x^\top  S^{-1} x + y^\top  S y.\label{young2}
\end{align}
\end{lemma}


\section{Observer design}
The aim of this section is to design a state observer for estimating the system state \eqref{main} and to obtain sufficient conditions for the existence of the proposed observer. The following observer is used to estimate the state of the system \eqref{main}
\begin{align}\label{observer}
\begin{split}
    z_{tt}(x,t) &=A_{1} z_{xx}(x,t)- B_{1} z_{t}(x,t)-D_{1} z(x,t)\\
    & \quad +  M f(\widehat{w}(x,t)) + G_{1} u(t) + L y(t),\\
    \widehat{w}(x,t)& =T z(x,t)+Q u(t)+E y(t),
\end{split}
\end{align}
subject to the Dirichlet boundary conditions $ z_{x}(0,t)=M w_{x}(0,t)$, $ z_{x}(1,t)=M w_{x}(1,t)$, and the initial conditions $(z_{0}(x),z_{1}(x))=(z(x,0), z_{t}(x,0))$ and $(\widehat{w}_{0}(x), \widehat{w}_{1}(x))=(\widehat{w}(x,0), \widehat{w}_{t}(x,0))$ where $z(x,t)\in [\mathcal{L}^{2}\big( (0,1); \mathbb R \big)]^{n}$, $\widehat{w}(x,t) \in [\mathcal{L}^{2}\big( (0,1); \mathbb R \big)]^{n}$ estimates $w(x,t)$, and $A_{1} \in \mathbb R^{n \times n}$, $B_{1} \in \mathbb R^{n \times n}$, $D_{1} \in \mathbb R^{n \times n}$,
$M \in \mathbb R^{n \times n}$,
$G_{1} \in \mathbb R^{n \times p}$, $L \in \mathbb R^{n \times q} $, $T \in \mathbb R^{n \times n} $, $Q \in \mathbb R^{n \times p}$, and $E \in \mathbb R^{n \times q} $ are unknown matrices and to be determined later such that the estimated error $e(x,t)=\widehat{w}(x,t)-w(x,t)$ converges to zero asymptotically. Now, motivated by \cite{delshad}, let us consider the error between $z(x,t)$ and $M w(x,t)$ as
\begin{equation}\label{er}
\varepsilon (x,t)=z(x,t)-M w(x,t).
\end{equation}
Therefore, we have
\begin{align}\label{er-01}
    \varepsilon _{tt}&(x,t)\cr
    =&z_{tt}(x,t)-M w_{tt}(x,t)\cr
    =&A_{1} z_{xx}(x,t)- B_{1} z_{t}(x,t)-D_{1} z(x,t)
   +M f(\widehat{w}(x,t))+G_{1} u(t)\cr
   & +L C \int_{0}^{1} w(x,t) dx + LK u(t)+LH  d(t)
   -MA w_{xx}(x,t)\cr
   & + MB w_{t}(x,t) + MD w(x,t)-M f(w(x,t))-MG u(t)-MF  d(t)\cr
   =&A_1 [z_{xx}(x,t)-M w_{xx}(x,t)]-B_1 [z_{t}(x,t)-M w_{t}(x,t)]\cr
   &  - D_1 [z(x,t)-M w(x,t)] + A_{1}M w_{xx}(x,t) - B_1M w_{t}(x,t) \cr & -D_1M w(x,t) + M f(\widehat{w}(x,t))+G_{1} u(t) +L C \int_{0}^{1} w(x,t) dx \cr
   & + LK u(t)+LH  d(t)
   -MA w_{xx}(x,t)  + MB w_{t}(x,t) \cr
   & + MD w(x,t)-M f(w(x,t))-MG u(t)-MF  d(t)\cr
   =&A_1 \varepsilon _{xx}(x,t)-B_1 \varepsilon _{t}(x,t)-D_1 \varepsilon (x,t)+[A_1M-MA]w_{xx}(x,t)\cr
   & - [B_1M-MB] w_{t}(x,t) - [D_1M-MD]w(x,t)+[LK+G_{1}-MG]u(t)\cr    %
   &+L C \int_{0}^{1} w(x,t) dx+M(f(\widehat{w} (x,t))-f(w(x,t))\cr
     &  +[LH-MF] d(t).
\end{align}
Also, the estimated error can be computed as
\begin{align}\label{er-11}
    e(x,t)& =\widehat{w}(x,t)-w(x,t)\cr
    &=Tz(x,t)+ Q u(t) + EC \int_{0}^{1} w(x,t) dx%
   + EK u(t)+ EH  d(t)- w(x,t)\cr
   &=T[z(x,t)-M w(x,t)]+ [Q+EK]u(t) +[TM-I_{n}]w(x,t)\cr 
      & \qquad \qquad +EH d(t) + EC \int_{0}^{1} w(x,t) dx \cr%
   &=T\varepsilon (x,t)+[TM-I_{n}]w(x,t) + [Q+EK]u(t) \cr
      & \qquad \qquad +EH d(t)+ EC \int_{0}^{1} w(x,t) dx.
\end{align}
Now, let $E=L$, so if there exist matrices $A_1$, $B_1$, $D_1$, $M$, $G_1$, $T$, $Q$, and $L$ such that
\begin{align}
\begin{cases}
  A_1M-MA=0,&\cr
  B_1M-MB=0,&\cr
  D_1M-MD=0,&\cr
  LH-MF=0, &\cr
  LK+G_{1}-MG=0,&\cr
 Q+L K=0,&\cr 
    LC=0,& \cr
 TM-I_{n}=0.& \label{matrix-eq}
 \end{cases}
\end{align}
 Hence, systems \eqref{er-01} and \eqref{er-11} get
\begin{align}
    \varepsilon _{tt}(x,t)
   &=A_1 \varepsilon _{xx}(x,t)-B_1 \varepsilon _{t}(x,t)-D_1 \varepsilon (x,t)+M \Delta f, \label{observer-error-1} \\ 
  e(x,t)
   &=T\varepsilon (x,t) +EH d(t), \label{observer-error-3}
\end{align}
 with $\Delta f=f(\widehat{w}(x,t))-f(w(x,t))$, the boundary conditions $ \varepsilon _{x}(0,t)=0$, and $ \varepsilon _{x}(1,t)=0$, and the initial conditions $\varepsilon_{0}(x)=\varepsilon(x,0)$, $\varepsilon_{1}(x)=\varepsilon_{t}(x,0)$. Also, equation \eqref{observer-error-3} leads $ e _{x}(0,t)=e _{x}(1,t)=0$, and the initial data of estimation error is $e_{0}(x)=T\varepsilon_{0}(x)+ EH d(0)$, and $e_{1}(x)=T\varepsilon_{1}(x)+ EH \dot{d}(0)$.
\section{Stability analysis}
For $d(t)=0$, the asymptotic stability of the estimation error $e(x,t)$ will be shown in the following theorem. Note that considering $d(t)=0$ leads to $e(x,t)=T \varepsilon(x,t)$; therefore, the asymptotic stability of $\varepsilon(x,t)$ follows the asymptotic stability of $e(x,t)$. 
\begin{theorem}\label{frst}
Suppose that the function $f$ is globally Lipschitz continuous on $[\mathcal{L}^{2}\big( (0,1); \mathbb R \big)]^{n}$
with 
Lipschitz constant $\gamma>0$,
 and let $d(t)=0$. Consider the system \eqref{main} together with the observer \eqref{observer}, and suppose that $A_{1}$, $B_{1}$, $D_{1}$, $M$, $P$, $T$ are positive definite matrices satisfying  the equations \eqref{matrix-eq} and the following inequalities
\begin{align}
&  \Pi= \begin{bmatrix}
    -B_{1}^\top P- P B_{1} + \delta P +\dfrac{\delta}{2} B_{1}^\top  B_{1} & PM  \\
    M^\top P & -\Gamma+\dfrac{\delta}{2}M^\top M 
    \end{bmatrix}\prec 0, \label{matrix-eq1}\\
& - \delta P D_{1}+ \delta P^{2} + \gamma^{2} T^\top  T\prec 0,~~ \Gamma -I \prec 0, \label{matrix-eq1-}
\end{align}    
    where $\delta < \min \{\frac{\lambda_{\min}(P)}{\lambda_{\max}(P)}, \frac{\lambda_{\min}(PD_{1})}{\lambda_{\max}(P)}\}$,
    and $P=P^\top $ is positive definite matrix, also $A_{1}$ and $D_{1}$ are symmetric matrices. Then, the observer error \eqref{observer-error-3} asymptotically approaches zero. 
\end{theorem}
\begin{proof}
Consider the following Lyapunov functional candidate
\begin{align}\label{lya}
 \Xi(t)= & \dfrac{1}{2}\int_{0}^{1}\varepsilon _{x}^\top (x,t) PA_1 \varepsilon _{x}(x,t)\, dx+\dfrac{1}{2}\int_{0}^{1}\varepsilon _{t}^\top (x,t) P \varepsilon _{t}(x,t)\, dx\cr
    &+\delta\int_{0}^{1}\varepsilon ^\top (x,t) P\varepsilon _{t}(x,t)\, dx+\dfrac{1}{2}\int_{0}^{1}\varepsilon ^\top (x,t) PD_1 \varepsilon (x,t)\, dx.   
\end{align}
The Lyapunov function $\Xi$ is nonnegative and satisfies the following inequality
\begin{align}
    C_{1}\| (\varepsilon,\varepsilon_{t})\|_{\mathcal{H}^{1}(0,1)\times \mathcal{L}^{2}(0,1)} \leq \Xi (t) \leq C_{2}\| (\varepsilon,\varepsilon_{t}) \|_{\mathcal{H}^{1}(0,1)\times \mathcal{L}^{2}(0,1)},
\end{align}
where
\begin{align*}
    C_{1}&=\min \{\frac{\lambda_{\min}(P A_{1})}{2},\frac{\lambda_{\min}(P)-\delta \lambda_{\max}(P)}{2},\frac{\lambda_{\min}(P D_{1})-\delta \lambda_{\max}(P)}{2}\},\\
    C_{2}&=\max \{\frac{\lambda_{\max}(P A_{1})}{2},\frac{\lambda_{\max}(P)+\delta \lambda_{\max}(P)}{2},\frac{\lambda_{\max}(P D_{1})+\delta \lambda_{\max}(P)}{2}\},
\end{align*}
with $\delta < \min \{\frac{\lambda_{\min}(P)}{\lambda_{\max}(P)},\frac{\lambda_{\min}(PD_{1})}{\lambda_{\max}(P)}\}$. Taking derivative respect to time yields
  \begin{align}\label{eq1}
 \dot{\Xi}(t)= &  \dfrac{1}{2}\int_{0}^{1}\varepsilon _{xt}^\top (x,t) PA_1 \varepsilon _{x}(x,t)\, dx+ \dfrac{1}{2}\int_{0}^{1}\varepsilon _{x}^\top (x,t) PA_1 \varepsilon _{xt}(x,t)\, dx\cr
 & +\dfrac{1}{2}\int_{0}^{1}\varepsilon _{tt}^\top (x,t) P \varepsilon _{t}(x,t)\, dx +\dfrac{1}{2}\int_{0}^{1}\varepsilon _{t}^\top (x,t) P \varepsilon _{tt}(x,t)\, dx\cr
    &+\delta\int_{0}^{1}\varepsilon _{t}^\top (x,t) P\varepsilon _{t}(x,t)\, dx +\delta\int_{0}^{1}\varepsilon ^\top (x,t) P\varepsilon _{tt}(x,t)\, dx\cr
    &+\dfrac{1}{2}\int_{0}^{1}\varepsilon _{t}^\top (x,t) PD_1 \varepsilon (x,t)\, dx +\dfrac{1}{2}\int_{0}^{1}\varepsilon ^\top (x,t) PD_1 \varepsilon _{t}(x,t)\, dx.    
\end{align}
Then, 
from \eqref{observer-error-1} it follows that
\begin{align}
 \dot{\Xi}(t)= &  \dfrac{1}{2}\int_{0}^{1}\varepsilon _{xt}^\top (x,t) PA_1 \varepsilon _{x}(x,t)\, dx+ \dfrac{1}{2}\int_{0}^{1}\varepsilon _{x}^\top (x,t) PA_1 \varepsilon _{xt}(x,t)\, dx\cr
 & +\dfrac{1}{2}\int_{0}^{1} \varepsilon _{xx}^\top (x,t)A_{1}^\top  P \varepsilon _{t}(x,t) \, dx -\dfrac{1}{2}\int_{0}^{1} \varepsilon _{t}^\top (x,t)B_{1}^\top  P \varepsilon _{t}(x,t) \, dx \cr
 &-\dfrac{1}{2}\int_{0}^{1} \varepsilon ^\top (x,t)D_{1}^\top  P \varepsilon _{t}(x,t) \, dx\
 + \dfrac{1}{2}\int_{0}^{1} \Delta f^\top  M^\top  P \varepsilon _{t}(x,t) \, dx\cr
 & +\dfrac{1}{2}\int_{0}^{1}\varepsilon _{t}^\top (x,t) PA_1 \varepsilon _{xx}(x,t)\, dx 
 -\dfrac{1}{2}\int_{0}^{1} \varepsilon _{t}^\top (x,t) P B_{1} \varepsilon _{t}(x,t) \, dx \cr
 &-\dfrac{1}{2}\int_{0}^{1} \varepsilon _{t}^\top (x,t) P D_{1} \varepsilon (x,t) \, dx + \dfrac{1}{2}\int_{0}^{1} \varepsilon _{t}^\top (x,t) P M \Delta f \, dx \cr
& +\delta\int_{0}^{1}\varepsilon _{t}^\top (x,t) P\varepsilon _{t}(x,t)\, dx 
+\delta\int_{0}^{1}\varepsilon ^\top (x,t) P A_{1}\varepsilon _{xx}(x,t)\, dx\cr
    &-\delta\int_{0}^{1}\varepsilon ^\top (x,t) P B_{1}\varepsilon _{t}(x,t)\, dx - \delta\int_{0}^{1}\varepsilon ^\top (x,t) P D_{1}\varepsilon (x,t)\, dx\cr
  & + \delta \int_{0}^{1} \varepsilon ^\top (x,t) P M \Delta f \, dx 
    +\dfrac{1}{2}\int_{0}^{1}\varepsilon _{t}^\top (x,t) PD_1 \varepsilon (x,t)\, dx \cr
    & +\dfrac{1}{2}\int_{0}^{1}\varepsilon ^\top (x,t) PD_1 \varepsilon _{t}(x,t)\, dx .
\end{align}
Under the assumptions on the matrices $P$ and $D_1$, we have $PD_{1}=D_{1}^\top P$, and
\begin{align}
 \dot{\Xi}(t)= &  \dfrac{1}{2}\int_{0}^{1}\varepsilon _{xt}^\top (x,t) PA_1 \varepsilon _{x}(x,t)\, dx+ \dfrac{1}{2}\int_{0}^{1}\varepsilon _{x}^\top (x,t) PA_1 \varepsilon _{xt}(x,t)\, dx\cr
 & +\dfrac{1}{2}\int_{0}^{1} \varepsilon _{xx}^\top (x,t)A_{1}^\top  P \varepsilon _{t}(x,t) \, dx\nonumber\\
 &-\dfrac{1}{2}\int_{0}^{1} \varepsilon _{t}^\top (x,t) [B_{1}^\top  P +P B_{1}  -\delta P] \varepsilon _{t}(x,t) \, dx \cr
& + \delta \int_{0}^{1} \varepsilon ^\top (x,t) P M \Delta f \, dx \cr
&
 + \dfrac{1}{2}\int_{0}^{1} \big[ \Delta f^\top  M^\top  P \varepsilon _{t}(x,t)+\varepsilon _{t}^\top (x,t) P M \Delta f \big] \, dx \cr
& +\dfrac{1}{2}\int_{0}^{1}\varepsilon _{t}^\top (x,t) PA_1 \varepsilon _{xx}(x,t)\, dx +\delta\int_{0}^{1}\varepsilon ^\top (x,t) P A_{1}\varepsilon _{xx}(x,t)\, dx\cr
    &-\delta\int_{0}^{1}\varepsilon ^\top (x,t) P B_{1}\varepsilon _{t}(x,t)\, dx - \delta\int_{0}^{1}\varepsilon ^\top (x,t) P D_{1}\varepsilon (x,t)\, dx. 
\end{align}
Now, using integration by parts, we get
\begin{align}
 \dot{\Xi}(t)= &  \dfrac{1}{2} \int_{0}^{1}\varepsilon _{xt}^\top (x,t) PA_1 \varepsilon _{x}(x,t)\, dx+ \dfrac{1}{2}\int_{0}^{1}\varepsilon _{x}^\top (x,t) PA_1 \varepsilon _{xt}(x,t)\, dx\cr
 &  + \dfrac{1}{2}  \varepsilon _{t}^\top (1,t) PA_1 \varepsilon _{x}(1,t)- \dfrac{1}{2} \varepsilon _{t}^\top (0,t) PA_1 \varepsilon _{x}(0,t)\cr
 &- \dfrac{1}{2}\int_{0}^{1}\varepsilon _{xt}^\top (x,t) PA_1 \varepsilon _{x}(x,t)\, dx\cr
 & + \dfrac{1}{2} \varepsilon _{x}^\top (1,t)A_{1}^\top  P \varepsilon _{t}(1,t) - \dfrac{1}{2} \varepsilon _{x}^\top (0,t)A_{1}^\top  P \varepsilon _{t}(0,t)\cr
& -\dfrac{1}{2}\int_{0}^{1} \varepsilon _{x}^\top (x,t)A_{1}^\top  P \varepsilon _{xt}(x,t) \, dx \cr
& + \delta \int_{0}^{1} \varepsilon ^\top (x,t) P M \Delta f \, dx \cr 
&  + \dfrac{1}{2}\int_{0}^{1} \big[ \Delta f^\top  M^\top  P \varepsilon _{t}(x,t)+\varepsilon _{t}^\top (x,t) P M \Delta f \big] \, dx \cr
& -\dfrac{1}{2}\int_{0}^{1} \varepsilon _{t}^\top (x,t) [ B_{1}^\top  P + P B_{1}-\delta P] \varepsilon _{t}(x,t) \, dx \cr
 &   -\delta\int_{0}^{1}\varepsilon ^\top (x,t) P B_{1}\varepsilon _{t}(x,t)\, dx\cr
    & - \delta\int_{0}^{1}\varepsilon ^\top (x,t) P D_{1}\varepsilon (x,t)\, dx
    + \delta \varepsilon ^\top (1,t) PA_1 \varepsilon _{x}(1,t)\cr
    & - \delta \varepsilon ^\top (0,t) PA_1 \varepsilon _{x}(0,t)
    - \delta \int_{0}^{1}\varepsilon _{x}^\top (x,t) P A_{1}\varepsilon _{x}(x,t)\, dx .
\end{align}
Next, using $P A_1=A_1^\top P$, we get from boundary conditions that
\begin{align}
 \dot{\Xi}(t)= 
&  \delta \int_{0}^{1} \varepsilon ^\top (x,t) P M \Delta f \, dx \cr
& + \dfrac{1}{2}\int_{0}^{1} \big[ \Delta f^\top  M^\top  P \varepsilon _{t}(x,t)+\varepsilon _{t}^\top (x,t) P M \Delta f \big] \, dx \cr
& -\dfrac{1}{2}\int_{0}^{1} \varepsilon _{t}^\top (x,t) [ B_{1}^\top  P + P B_{1}-\delta P] \varepsilon _{t}(x,t) \, dx \cr
  &  -\delta\int_{0}^{1}\varepsilon ^\top (x,t) P B_{1}\varepsilon _{t}(x,t)\, dx\cr
    & - \delta \int_{0}^{1}\varepsilon ^\top (x,t) P D_{1}\varepsilon (x,t)\, dx
    - \delta \int_{0}^{1}\varepsilon _{x}^\top (x,t) P A_{1}\varepsilon _{x}(x,t)\, dx .
\end{align}
It can be derived from inequalities \eqref{young1}-\eqref{young2}, with $S=I$, that
\begin{align}
      \delta &\int_{0}^{1} \varepsilon ^\top (x,t) P M \Delta f \, dx \cr
       & \leq\dfrac{\delta}{2} \int_{0}^{1} \varepsilon ^\top (x,t) P^{2} \varepsilon (x,t) \, dx + \dfrac{\delta}{2}  \int_{0}^{1} \Delta f^\top  M^\top  M \Delta f \, dx, \label{y-1}\\
      -\delta &\int_{0}^{1}\varepsilon ^\top (x,t) P B_{1}\varepsilon _{t}(x,t)\, dx  \cr
       & \leq \dfrac{\delta}{2}\int_{0}^{1}\varepsilon ^\top (x,t) P^{2} \varepsilon (x,t)\, dx + \dfrac{\delta}{2} \int_{0}^{1}\varepsilon _{t}^\top (x,t) B_{1}^\top  B_{1}\varepsilon _{t}(x,t)\, dx. \label{y-2}
\end{align}
 Using \eqref{y-1} and \eqref{y-2}, we have
\begin{align}
    \dot{\Xi}(t) \leq 
 & - \delta \int_{0}^{1} \varepsilon ^\top (x,t) [P D_{1}-P^{2}]\varepsilon (x,t) \, dx \cr
& + \dfrac{1}{2}\int_{0}^{1} \big[ \Delta f^\top  M^\top  P \varepsilon _{t}(x,t)+\varepsilon _{t}^\top (x,t) P M \Delta f + \dfrac{\delta}{2}\Delta f^\top  M^\top  M \Delta f \big] \, dx \cr
& -\dfrac{1}{2}\int_{0}^{1} \varepsilon _{t}^\top (x,t) [ B_{1}^\top  P + P B_{1}-\delta P -\dfrac{\delta}{2} B_{1}^\top  B_{1}] \varepsilon _{t}(x,t) \, dx  \cr
& - \delta \int_{0}^{1}\varepsilon _{x}^\top (x,t) P A_{1}\varepsilon _{x}(x,t)\, dx .
\end{align}
Thanks to the global Lipschitz continuity of $f(w(x,t))$, and since $\Gamma \prec I$, we get
\begin{align}\label{eq-lip}
\Delta f^\top  \Gamma \Delta f & \leq \gamma^{2} (w(x,t) - \widehat{w} (x,t))^\top   (w(x,t) - \widehat{w} (x,t))\cr
& = \gamma^{2} e^\top (x,t)  e(x,t) =\gamma^{2} \varepsilon ^\top (x,t) T^\top T \varepsilon(x,t).
\end{align}
Adding and subtracting the term $\int_{0}^{1} \Delta f^\top  \Gamma \Delta f dx$ yields
\begin{align}
    \dot{\Xi}(t) \leq 
 & - \delta \int_{0}^{1} \varepsilon ^\top (x,t) [P D_{1}-P^{2}]\varepsilon (x,t) \, dx + \int_{0}^{1} [\Delta f^\top  \Gamma \Delta f - \Delta f^\top  \Gamma \Delta f]\, dx \cr
& + \dfrac{1}{2}\int_{0}^{1} \big[ \Delta f^\top  M^\top  P \varepsilon _{t}(x,t)+\varepsilon _{t}^\top (x,t) P M \Delta f + \dfrac{\delta}{2}\Delta f^\top  M^\top  M \Delta f \big] \, dx \cr
& -\dfrac{1}{2}\int_{0}^{1} \varepsilon _{t}^\top (x,t) [ B_{1}^\top  P + P B_{1}-\delta P -\dfrac{\delta}{2} B_{1}^\top  B_{1}] \varepsilon _{t}(x,t) \, dx  \cr
& - \delta \int_{0}^{1}\varepsilon _{x}^\top (x,t) P A_{1}\varepsilon _{x}(x,t)\, dx.
\end{align}
Thus, considering \eqref{eq-lip} leads to
\begin{align}
    \dot{\Xi}(t) \leq 
 & -  \int_{0}^{1} \varepsilon ^\top (x,t) [\delta P D_{1}- \delta P^{2}-\gamma^{2} T^\top  T]\varepsilon (x,t) \, dx \cr
&  + \dfrac{1}{2}\int_{0}^{1} \big[  \Delta f^\top  M^\top  P \varepsilon _{t}(x,t) + \varepsilon _{t}^\top (x,t) P M \Delta f + \Delta f^\top  (-\Gamma\nonumber \\
&+ \dfrac{\delta}{2} M^\top  M) \Delta f \big] \, dx \cr
& -\dfrac{1}{2}\int_{0}^{1} \varepsilon _{t}^\top (x,t) [ B_{1}^\top  P + P B_{1}-\delta P -\dfrac{\delta}{2} B_{1}^\top  B_{1}] \varepsilon _{t}(x,t) \, dx \nonumber\\
 &- \delta  \int_{0}^{1}\varepsilon _{x}^\top (x,t)  P A_{1}\varepsilon _{x}(x,t)\, dx\cr
= &- \delta \int_{0}^{1} \varepsilon ^\top (x,t) [\delta P D_{1}- \delta P^{2}-\gamma^{2} T^\top  T]\varepsilon (x,t) \, dx\nonumber \\
&+ \dfrac{1}{2}\int_{0}^{1} [ \varepsilon _{t}^\top (x,t) \, \Delta f^\top  ] \Pi \begin{bmatrix}
     \varepsilon _{t}(x,t) \\ \Delta f 
\end{bmatrix}\, dx  \cr
& - \delta \int_{0}^{1}\varepsilon _{x}^\top (x,t) P A_{1} \varepsilon _{x}(x,t)\, dx.
\end{align}
Therefore, if the inequalities \eqref{matrix-eq1}--\eqref{matrix-eq1-} hold, then the time derivative of Lyapunov function $\Xi(t)$ is negative definite, which shows the state estimation error goes to zero when time goes to infinity and completes the proof.
\end{proof}
\section{$H_{\infty}$ design}
In this section, we solve the following  $H_{\infty}$ observer design problem: given a prescribed attenuation level $\mu >0$, design an observer in a way that the following conditions are fulfilled:
\begin{itemize}
\item[(C1)] The estimation error state \eqref{observer-error-3} with $d(t)=0$ is asymptotically stable.

\item[(C2)] Under zero initial condition, the inequality 
\begin{equation}
\int_{0}^{+\infty} e^\top (x,t) e(x,t) \,dt \leq \mu^{2} \int_{0}^{+\infty} d^\top (t) d(t) \,dt
\end{equation}
 holds with the non-zero disturbance input $d(t)$.
 \end{itemize}
\begin{theorem}
Let $f$ be globally Lipschitz continuous on $[\mathcal{L}^{2}\big( (0,1); \mathbb R \big)]^{n}$. Given a scalar $\mu >0$, then the estimation error state \eqref{observer-error-3} produced by the observer \eqref{observer} for positive Lipschitz constant $\gamma$ is asymptotically stable with $H_{\infty}$ performance level $\mu>0$, if there exist positive definite matrices $A_{1}$, $B_{1}$, $D_{1}$, $M$, $P$, $T$ that satisfy in the equations \eqref{matrix-eq}, and the following inequalities
\begin{align}\label{matrix-eq2}
& \Theta=\begin{bmatrix}
  a_{11} & 0 & PM & 0  \\
   0 & a_{22} & 0 & (1+\gamma^{2})T^\top EH\\
   M^\top P & 0 & a_{33} & 0\\
    0 & (1+\gamma^{2})H^\top E^\top T & 0 & a_{44}
    \end{bmatrix} \prec 0,
~~\Gamma - I \prec 0,
\end{align}    
with
\begin{align*}
a_{11}&= - B_{1}^\top P- P B_{1} + \delta P + \dfrac{\delta}{2} B_{1}^\top  B_1,\\
a_{22}&= -\delta PD_{1} +\delta P^{2} + (1+\gamma^{2})T^\top T,\\
a_{33}&=-\Gamma +\dfrac{\delta}{2} M^\top M,\\
a_{44}&=(1+\gamma^{2})H^\top E^\top EH -\mu^{2} I,
\end{align*}
    where $\delta < \min  \{ \frac{\lambda_{\min}(P)}{\lambda_{\max}(P)}, \frac{ \lambda_{\min}(PD_{1})}{\lambda_{\max}(P)}\}$,
    and $P$, $A_{1}$ and $D_{1}$ are symmetric matrices. 
\end{theorem}
\begin{proof}
The asymptotic stability of \eqref{observer-error-3} when $d(t)=0$ can be achieved by Theorem \ref{frst}. Now, for the non-zero input $d(t)$, we will show under zero initial condition, the inequality $\| e(x,t)\|_{\mathcal{L}^{2}(0,\infty)}^{2} \leq \mu^{2} \| d(t)\|_{\mathcal{L}^{2}(0,\infty)}^{2}$ holds. 
Consider the Lyapunov function $\Xi(t)$ \eqref{lya}, then using the similar arguments to those used for equations \eqref{eq1}–\eqref{y-2}, it can be follows that
\begin{align}
    \dot{\Xi}(t) \leq 
 & - \delta \int_{0}^{1} \varepsilon ^\top (x,t) [P D_{1}-P^{2}]\varepsilon (x,t) \, dx \cr
& + \dfrac{1}{2}\int_{0}^{1} \big[ \Delta f^\top  M^\top  P \varepsilon _{t}(x,t)+\varepsilon _{t}^\top (x,t) P M \Delta f + \dfrac{\delta}{2}\Delta f^\top  M^\top  M \Delta f \big] \, dx \cr
& -\dfrac{1}{2}\int_{0}^{1} \varepsilon _{t}^\top (x,t) [ B_{1}^\top  P + P B_{1}-\delta P -\dfrac{\delta}{2} B_{1}^\top  B_{1}] \varepsilon _{t}(x,t) \, dx  \cr
& - \delta \int_{0}^{1}\varepsilon _{x}^\top (x,t) P A_{1}\varepsilon _{x}(x,t)\, dx\cr
& = - \delta \int_{0}^{1} \varepsilon ^\top (x,t) [P D_{1}-P^{2}]\varepsilon (x,t) \, dx + \int_{0}^{1} [\Delta f^\top  \Gamma \Delta f - \Delta f^\top  \Gamma \Delta f]\, dx \cr
& + \dfrac{1}{2}\int_{0}^{1} \big[ \Delta f^\top  M^\top  P \varepsilon _{t}(x,t)+\varepsilon _{t}^\top (x,t) P M \Delta f + \dfrac{\delta}{2}\Delta f^\top  M^\top  M \Delta f \big] \, dx \cr
& -\dfrac{1}{2}\int_{0}^{1} \varepsilon _{t}^\top (x,t) [ B_{1}^\top  P + P B_{1}-\delta P -\dfrac{\delta}{2} B_{1}^\top  B_{1}] \varepsilon _{t}(x,t) \, dx  \cr
& - \delta \int_{0}^{1}\varepsilon _{x}^\top (x,t) P A_{1}\varepsilon _{x}(x,t)\, dx .
\end{align}
 It can be derived from \eqref{observer-error-3} that
 \begin{align}
 e^\top (x,t)  e(x,t) & = \varepsilon ^\top (x,t) T^\top T \varepsilon (x,t) +\varepsilon ^\top (x,t)T^\top EH d \cr
& +d^\top H^\top E^\top T \varepsilon(x,t)+ d^\top H^\top E^\top EH d .
 \end{align}
Furthermore, due to the global Lipschitz of $f(w(x,t))$, we get
\begin{align}\label{eq-lip-1}
\Delta f^\top   \Gamma \Delta f & \leq \gamma^{2} (\widehat{w}(x,t)-w(x,t))^\top   (\widehat{w}(x,t)-w(x,t)) \cr
 & = \gamma^{2} \varepsilon ^\top (x,t) T^\top T \varepsilon (x,t) + \gamma^{2} \varepsilon ^\top (x,t)T^\top EH d\cr
& + \gamma^{2} d^\top H^\top E^\top T\varepsilon(x,t)+ \gamma^{2} d^\top H^\top E^\top EH d. 
\end{align}
Therefore, we get
\begin{align}
    \dot{\Xi}(t) & +\int_{0}^{1}[e^\top (x,t)e(x,t)-\mu^{2} d^\top  d]\, dx \cr
  \leq & -  \int_{0}^{1} \varepsilon ^\top (x,t) [\delta P D_{1}- \delta P^{2}- (1+\gamma^{2}) T^\top  T]\varepsilon (x,t) \, dx  -\int_{0}^{1} \mu^{2} d^\top  d \, dx \cr
& + \dfrac{1}{2}\int_{0}^{1} \big[\Delta f^\top  M^\top  P \varepsilon _{t}(x,t) \cr
& \qquad\qquad + \varepsilon _{t}^\top (x,t) P M \Delta f + \Delta f^\top  (- \Gamma + \dfrac{\delta}{2} M^\top  M) \Delta f \big] \, dx \cr
&+ \int_{0}^{1} [(1+\gamma^{2})\varepsilon ^\top (x,t)T^\top EH d + (1+\gamma^{2}) d^\top H^\top E^\top T\varepsilon(x,t)  \cr
& \qquad\qquad+ (1+\gamma^{2}) d^\top H^\top E^\top EH d]\, dx  \cr
& -\dfrac{1}{2}\int_{0}^{1} \varepsilon _{t}^\top (x,t) [ B_{1}^\top  P + P B_{1}-\delta P -\dfrac{\delta}{2} B_{1}^\top  B_{1}] \varepsilon _{t}(x,t) \, dx  \cr
& - \delta \int_{0}^{1}\varepsilon _{x}^\top (x,t)  P A_{1}\varepsilon _{x}(x,t)\, dx\cr
& = -\delta \int_{0}^{1} \varepsilon_{x} ^\top (x,t)  P A_{1} \varepsilon_{x} (x,t) \, dx 
+ \dfrac{1}{2}\int_{0}^{1} \mathcal{Y}^\top (x,t)  \Theta \mathcal{Y}(x,t)\, dx 
\end{align}
with $\mathcal{Y}(x,t)=\begin{bmatrix}
\varepsilon _{t}^\top (x,t) & \varepsilon ^\top (x,t) &  \Delta f^\top  & d^\top 
\end{bmatrix}^\top $.

Therefore, if $\Theta\prec 0$, then thanks to $P A_{1} \succ 0$, we have
\begin{equation}\label{lya-eq}
 \dot{\Xi}\leq \int_{0}^{1}[\mu^{2} d^\top  d - e^\top (x,t)e(x,t)]\, dx.   
\end{equation}
Now, taking integral both side of \eqref{lya-eq} from $0$ to $\infty$ respect to $t$ leads to
\begin{equation}
\Xi(x,\infty)-\Xi(x,0)\leq \int_{0}^{1}[\mu^{2} \| d(t)\|_{\mathcal{L}^{2}(0,\infty)}^{2} -\| e(x,t) \|_{\mathcal{L}^{2}(0,\infty)}^{2}]\, dx, 
\end{equation}
so, for zero initial condition, we obtain
\begin{equation}
\Xi(x,\infty) \leq \int_{0}^{1}[\mu^{2} \| d(t)\| _{\mathcal{L}^{2}(0,\infty)}^{2} -\| e(x,t) \| _{\mathcal{L}^{2}(0,\infty)}^{2} ]\, dx, 
\end{equation}
which leads to $\| e(x,t) \| _{\mathcal{L}^{2}(0,\infty)}^{2} \leq \mu^{2} \| d(t)\|_{\mathcal{L}^{2}(0,\infty)}^{2}$ that shows estimation error state is asymptotically stable with $H_{\infty}$ performance level $\mu$ and completes the proof. 
\end{proof}
\section{Numerical simulation}
  \begin{figure}
\subfloat[$w_{1}(x,t)$] {\includegraphics[scale=.41]{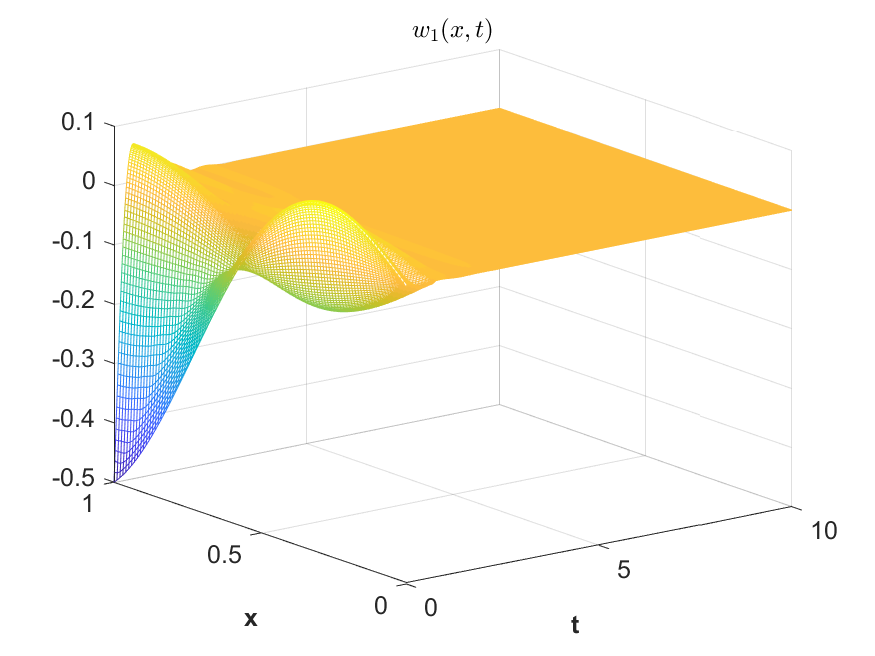}}
\subfloat[$w_{2}(x,t)$]{\includegraphics[scale=.41]{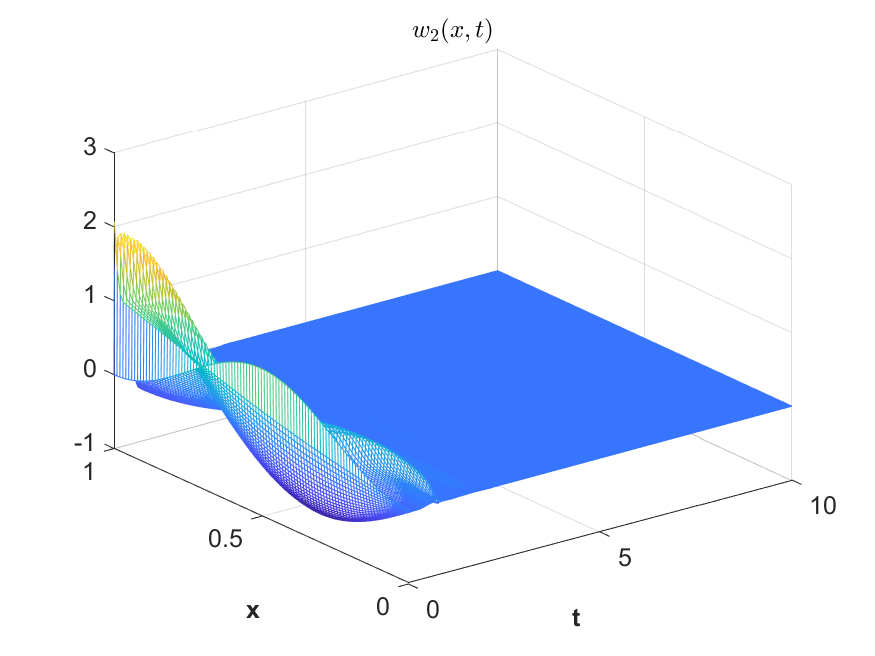}}\\
\subfloat[$\widehat{w}_{1}(x,t)$] {\includegraphics[scale=.41]{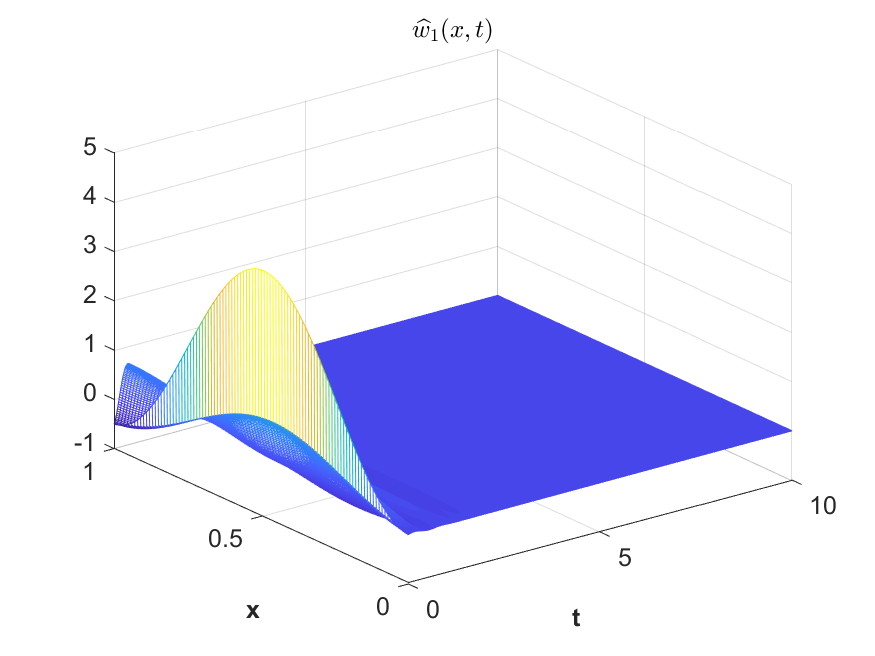}}
\subfloat[$\widehat{w}_{2}(x,t)$]{\includegraphics[scale=.41]{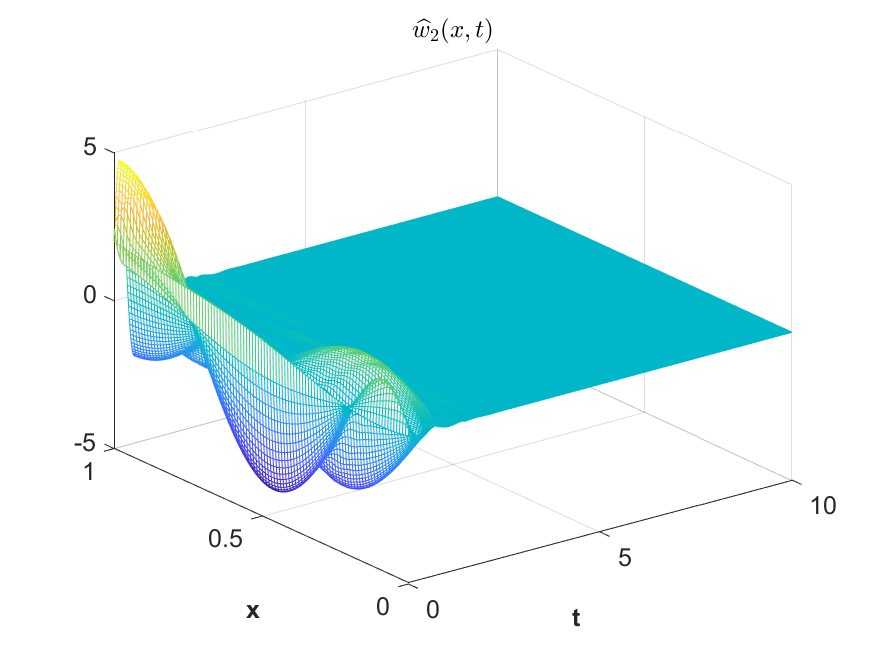}}\\
\subfloat[$e_1(x,t)$]{\includegraphics[scale=.41]{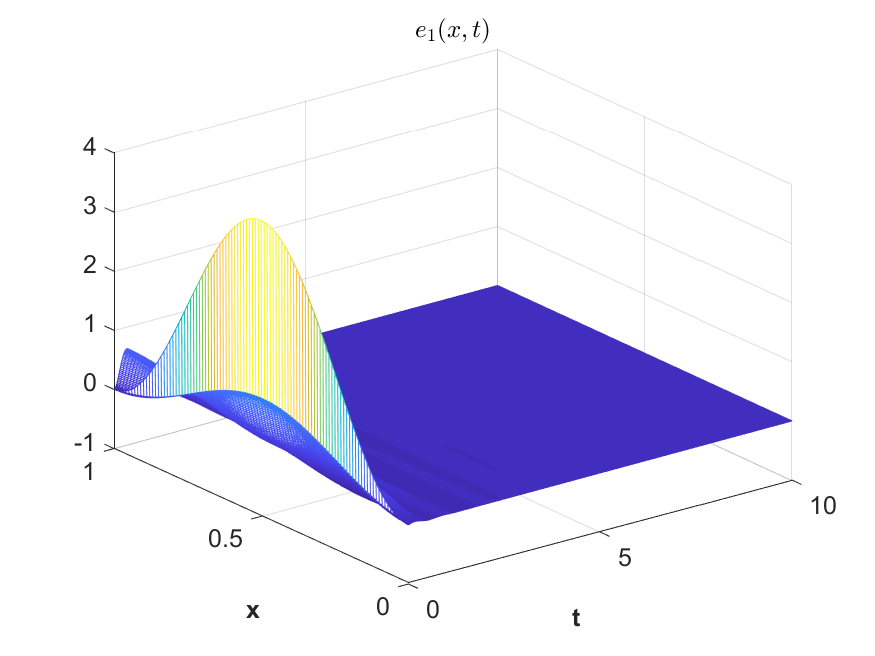}}
\subfloat[$e_2(x,t)$]{\includegraphics[scale=.41]{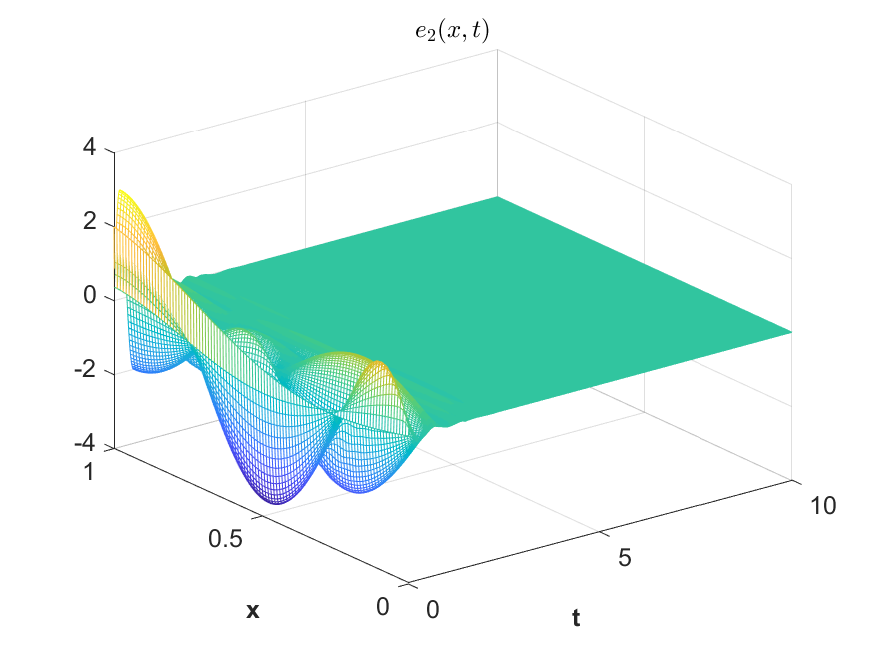}}
\caption{Response of the original open-loop plant \eqref{main}, the estimated system \eqref{observer}, the error $e(x,t)=\widehat{w}(x,t)-w(x,t)$ with no disturbance ($d(t)=0$).}\label{figfn}
\end{figure} 
 \begin{figure}
\subfloat[$w_{1}(x,t)$] {\includegraphics[scale=.41]{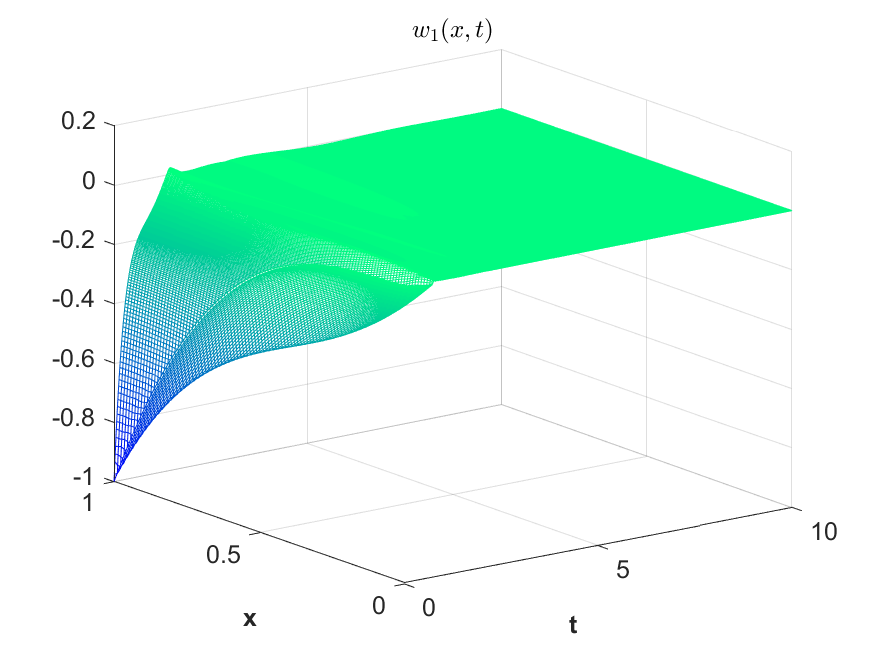}}
\subfloat[$w_{2}(x,t)$]{\includegraphics[scale=.41]{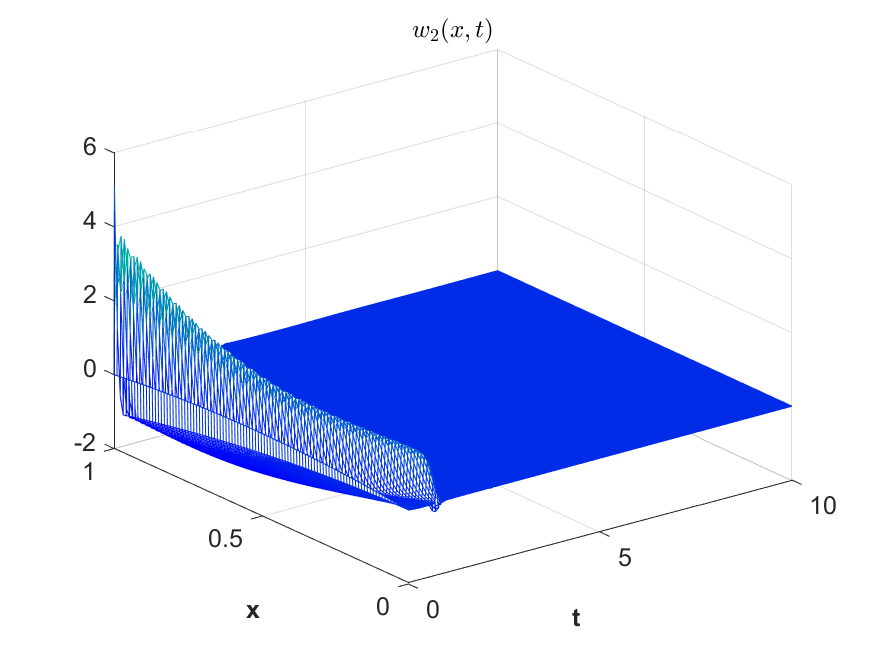}}\\
\subfloat[$\widehat{w}_{1}(x,t)$] {\includegraphics[scale=.41]{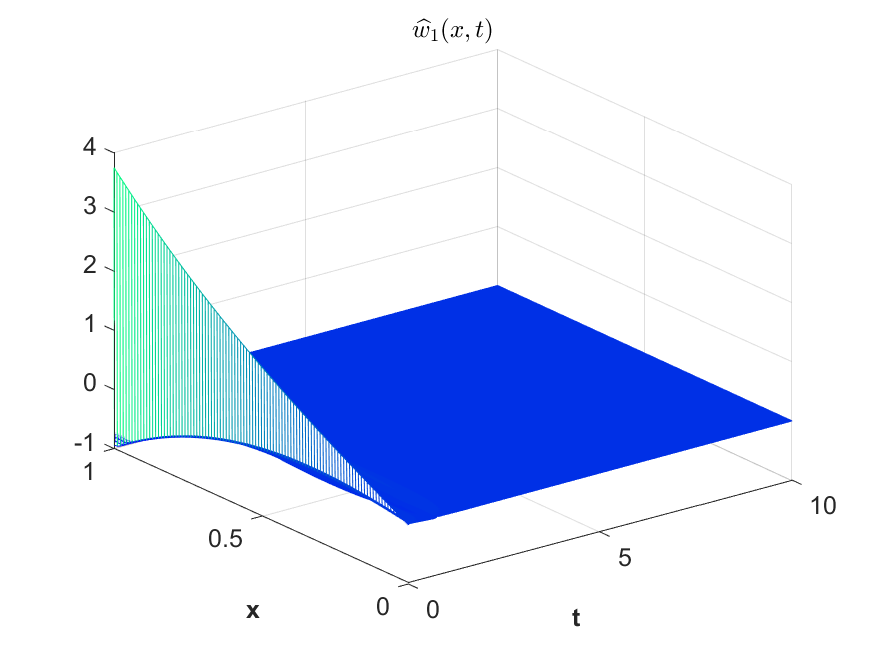}}
\subfloat[$\widehat{w}_{2}(x,t)$]{\includegraphics[scale=.41]{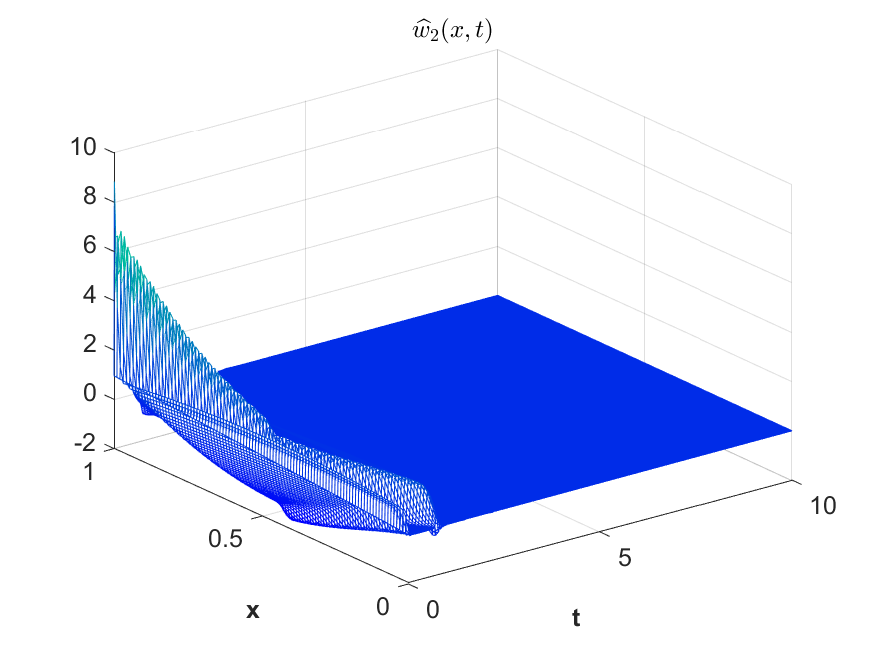}}\\
\subfloat[$e_1(x,t)$]{\includegraphics[scale=.41]{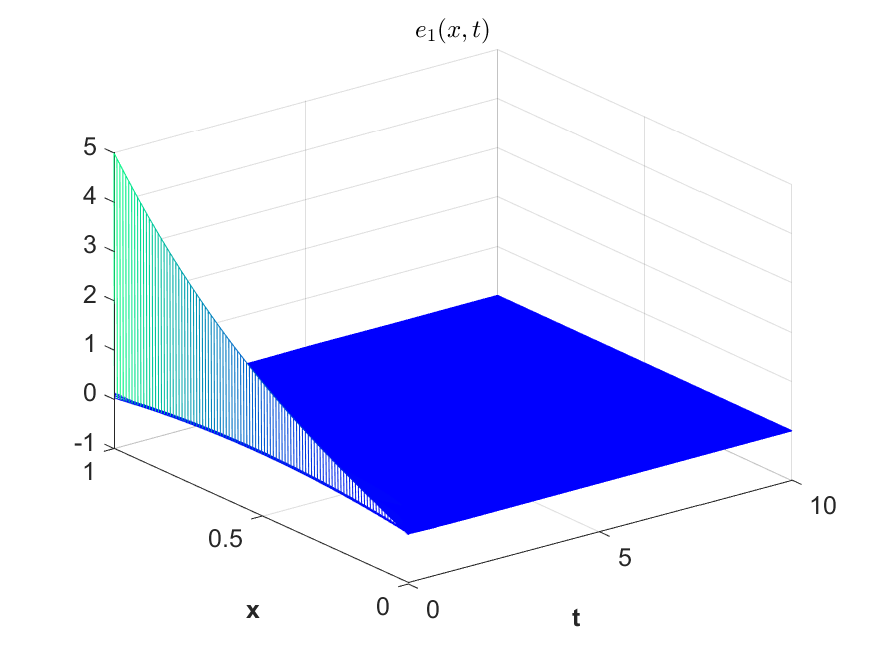}}
\subfloat[$e_2(x,t)$]{\includegraphics[scale=.41]{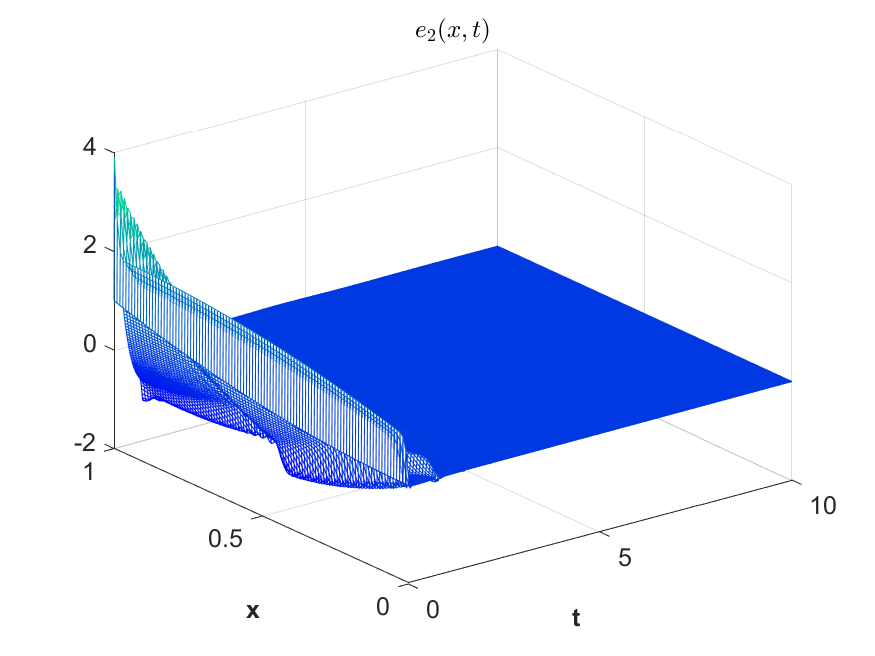}}
\caption{Response of the original open-loop plant \eqref{main}, the estimated system \eqref{observer}, the error system $e(x,t)=\widehat{w}(x,t)-w(x,t)$ in the presence of disturbance ($d(t) \neq 0$).}\label{figsn}%
\end{figure} 
In this section, we present some numerical simulations to illustrate the effectiveness of the proposed method. Consider the system \eqref{main} of two ($n = 2$) coupled wave equations
with, $q=2$, $p=d=1$, and
\begin{align}
  &  A=\begin{bmatrix}
     2 & 0\\
     0 & 2
    \end{bmatrix}, ~~ B=\begin{bmatrix}
     4 & 0\\
     0 & 4
    \end{bmatrix},~~D=\begin{bmatrix}
     4.5 & 0\\
     0 & 4.5
    \end{bmatrix},\cr
   & G=\begin{bmatrix}
     2  \\
     2  
\end{bmatrix},~
C=\begin{bmatrix}
     0.75 & -0.75\\
     -0.75 & 0.75
    \end{bmatrix},~K=\begin{bmatrix}
    - 2  \\
     -2  
    \end{bmatrix}\cr
    & H=\begin{bmatrix}
     0.05 \\
     0.05
    \end{bmatrix},
    ~ F=\begin{bmatrix}
     0.1   \\
     -0.1   
    \end{bmatrix},~ C_{1}=\begin{bmatrix}
     0 & 0.25\\
     0.25 & 0
    \end{bmatrix}
\end{align}
Also, let $f(w(x,t))= \begin{bmatrix}
0.1 \sin(w_{1}(x,t)) & 0.1 \sin(w_{2}(x,t)) 
\end{bmatrix}^\top $.
For simulation purpose, we assume $u(t)= - K_{1} \int_{0}^{1} w(x,t)\, dx$ with  $K_1=\begin{bmatrix}
     \frac{1}{2}& \frac{1}{2}
    \end{bmatrix}$, and consider the problem in the presence of disturbance (i.e., $d(t)\neq 0$) and without it. Solving the matrix inequalities \eqref{matrix-eq1} and \eqref{matrix-eq2} respect to constraints of \eqref{matrix-eq} yields 
\begin{align}
  &  A_{1}=\begin{bmatrix}
     2 & 0\\
     0 & 2
    \end{bmatrix}, ~~ B_{1}=\begin{bmatrix}
     4 & 0\\
     0 & 4
    \end{bmatrix},~~D_{1}=\begin{bmatrix}
     4.5 & 0\\
     0 & 4.5
    \end{bmatrix},\cr
   & G_{1}=\begin{bmatrix}
     5.7 \\
     - 1.9   
\end{bmatrix},~
M=\begin{bmatrix}
     0.95 & 0\\
     0 & 0.95
    \end{bmatrix},~T=\begin{bmatrix}
     100/95 & 0\\
     0 & 100/95
    \end{bmatrix}\cr
    & L=\begin{bmatrix}
     0.95 & 0.95\\
     -0.95 & -0.95
    \end{bmatrix},~ Q=\begin{bmatrix}
     3.8  \\
     -3.8  
    \end{bmatrix}, 
   E=\begin{bmatrix}
     0.95 & 0.95\\
     -0.95 & -0.95
    \end{bmatrix},\cr
    & P=\begin{bmatrix}
     2.5 & 0\\
     0 & 2.5
    \end{bmatrix},~ \Gamma=\begin{bmatrix}
     8/9 & 0\\
     0 & 8/9
    \end{bmatrix}
\end{align}
also, $\gamma=0.1$, $\delta=1/4$, $\mu=0.95$. 

To solve the considered PDE systems, the semi-discretization method is
used in all simulations by discretizing the spatial solution domain $x \in [0,1]$ and solving the resulting discretized system of ODEs in the Matlab environment by using the Runge-Kutta method with $dx=dt=0.01$.
When $d(t)=0$, the initial conditions are set as
\begin{align*}
& w_{0}(x) =\begin{bmatrix}
  -0.5 \cos( 2\pi x) +0.5 & 0.5 x \cos(\pi x)
\end{bmatrix}^\top ,~ w_{1}(x) =\begin{bmatrix}
 x^2 - x  & x
\end{bmatrix}^\top ,\\
& \widehat{w}_{0}(x) =\begin{bmatrix}
 -2.5 \cos( 2\pi x) +2.5 & - 1.5 x \cos(\pi x)
\end{bmatrix}^\top ,~ \widehat{w}_{1}(x) =\begin{bmatrix}
 -x^2 - x & - x
\end{bmatrix}^\top ,\\
& e_{0}(x) =\begin{bmatrix}
 -2 \cos( 2\pi x) +2 & -2 x \cos(\pi x)
\end{bmatrix}^\top ,~ e_{1}(x) =\begin{bmatrix}
 - 2 x^2 & -2 x
\end{bmatrix}^\top ,
\end{align*}
and for $d(t)=0.2 e^{-0.4 t} \sin(0.5 \pi t)$, the initial conditions are set as
\begin{align*}
& w_{0}(x) =\begin{bmatrix}
 x-x^2 & -x^3
\end{bmatrix}^\top ,~ w_{1}(x) =\begin{bmatrix}
- x^2+2x-1 & -x^2
\end{bmatrix}^\top ,\\
& \widehat{w}_{0}(x) =\begin{bmatrix}
4 x +2  x^3 -x^2 &  0.25 x^3 - x^2
\end{bmatrix}^\top ,~
 \widehat{w}_{1}(x) =\begin{bmatrix}
0 &  -x3 + x^2
\end{bmatrix}^\top ,\\
& e_{0}(x) =\begin{bmatrix}
3 x +2 x^3 & x^2
\end{bmatrix}^\top ,~ e_{1}(x) =\begin{bmatrix}
  x^2-2x+1 & 0.25 x^3
\end{bmatrix}^\top .
\end{align*}
In absence of disturbance, Figure \eqref{figfn} exhibits the solution of original system \eqref{main}, the estimated state \eqref{observer} and error system $e(x,t)$ (i.e., equation \eqref{observer-error-3})  where the simulations for $w_{1}(x,t)$, $w_{2}(x,t)$, $\widehat{w}_1(x,t)$, $\widehat{w}_1(x,t)$, $e_{1}(x,t)$ and $e_{2}(x,t)$ are plotted in Figures \ref{figfn}(a), \ref{figfn}(b), \ref{figfn}(c), \ref{figfn}(d), \ref{figfn}(e) and \ref{figfn}(f), respectively. Furthermore, the solution of original system \eqref{main}, the estimated state \eqref{observer} and the error of between them in the presence of the disturbance are depicted in Figure \ref{figsn} in which the simulations for $w_{1}(x,t)$, $w_{2}(x,t)$, $\widehat{w}_1(x,t)$, $\widehat{w}_1(x,t)$, $e_{1}(x,t)$ and $e_{2}(x,t)$ are plotted in Figures \ref{figsn}(a), \ref{figsn}(b), \ref{figsn}(c), \ref{figsn}(d), \ref{figsn}(e) and \ref{figsn}(f), respectively. When the disturbance occurs, the system behaviour is shown in Figure \ref{figsn}, which shows that the $H_{\infty}$ performance in the presence of disturbance is satisfied. These figures confirm that the proposed observer estimates the state quite efficiently for uncertain wave systems with unknown input.

\section{Conclusion}
In this paper, we introduce an $H_{\infty}$ observer designed for a class of coupled uncertain Lipschitz nonlinear wave PDE systems. The observer structure and the dynamic error components, which include undetermined parameters, are constructed for these wave equations under the condition $ d(t) = 0$ and for the case $d(t)\neq 0$, with $ \sup_{0 \neq d(t) \in \mathcal{L}^{2}(0,+\infty)} \dfrac{\| e(x, t) \|_{\mathcal{L}^{2}}}{\| d(t) \|_{\mathcal{L}^{2}}} < \mu$. The sufficient conditions for the proposed method, ensuring asymptotic stability, are provided, and the $H_{\infty}$ observer performance is examined, with the stability analysis and necessary conditions expressed through matrix inequalities. A numerical example demonstrates the efficacy of the proposed method, although the conditions specified in \eqref{matrix-eq} are quite restrictive. 

\section*{Acknowledgement} This work was carried out at the University of Wuppertal and was funded by the German Academic Exchange Service (DAAD) under grant number 57681230. 
\section*{Conflict of Interest} The authors declare no potential conflicts of interest.

\section*{Data availability statement}
Data sharing is not applicable to this article as no new data were created or analyzed in this study.


\end{document}